\documentstyle[12pt]{amsart}
\def\cases#1{\left\{\,\vcenter{\normalbaselines
    \ialign{$##\hfil$&\quad##\hfil\crcr#1\crcr}}\right.}

\newtheorem{theorem}{THEOREM}[section]

\newtheorem{corollary}[theorem]{Corollary}
\newtheorem{lemma}[theorem]{Lemma}

\newtheorem{proposition}[theorem]{Proposition}
\def\hbar{\overline{h}}

\def\CC{{\rm\kern.24em\vrule
width.02em height1.4ex
depth-.05ex\kern-.26em C}}
\def\QQ{{\rm\kern.24em\vrule width.02em
height1.4ex depth-.05ex\kern-.26em Q}}
\def\RR{{\rm I\kern-.2em R}}
\def\HH{{\rm I\kern-.2em H}}
\def\ZZ{{\rm\kern.26em\vrule width.02em
height0.5ex depth0ex\kern.04em\vrule width.02em
height1.47ex depth-1ex\kern-.34em Z}}
\def\Ibb#1{{\rm I\kern-.23em#1}}
\def\Ib#1{{\rm I\kern-.25em#1}}
\def\k#1{\kern#1em}
\def\vb#1{\vrule width.02em height1.4ex depth-.05ex}

\def\HH{Inn H}
\def\11{{\rm\k{.45}\vb0\k{-.142}1}}

\def\epf{\hskip.2in\vrule width.4pt height6.65pt
depth.15pt\vrule
width2.5pt height6.65pt depth-6.25pt\hskip-2.5pt\vrule
width2.5pt
height.25pt depth.15pt\vrule width.4pt
height6.65pt depth.15pt\ }

\def\proof{\noindent {\bf Proof. }}

\def \d{\partial}

\def \hbar{\overline{h}}

\def\supp{\hbox{supp}}

\font\teneufm=eufm10
\font\seveneufm=eufm7
\font\fiveeufm=eufm5
\newfam\eufmfam
\textfont\eufmfam=\teneufm
\scriptfont\eufmfam=\seveneufm
\scriptscriptfont\eufmfam=\fiveeufm

\newfam\msbfam
\font\tenmsb=msbm10 scaled \magstep1   \textfont\msbfam=\tenmsb
\font\sevenmsb=msbm7 scaled \magstep1    \scriptfont\msbfam=\sevenmsb
\font\fivemsb=msbm5 scaled \magstep1     \scriptscriptfont\msbfam=\fivemsb
\def\Bbb{\fam\msbfam \tenmsb}

\def\RR{{\Bbb R}}
\def\CC{{\Bbb C}}
\def\QQ{{\Bbb Q}}

\def\ZZ{{\Bbb Z}}
\def\HH{{\Bbb H}}

\def\Dbar{\overline{D}}

\begin{document}

\title[Boundary Value Problems
for the Inhomogeneous Laplace Equation]
{Elliptic Boundary Value Problems\\
for the Inhomogeneous Laplace Equation on Bounded Domains}

\author{Steven G. Krantz}
\author{Song-Ying Li}
\address{\hskip-\parindent
Steven G. Krantz\\
Department of Mathematics\\
Washington University\\
St. Louis, MO 63130}
\email{sk@@math.wustl.edu}
\address{\hskip-\parindent
Song-Ying Li\\
Department of Mathematics\\
Washington University\\
St. Louis, MO 63130}
\email{songying@@math.wustl.edu}

\thanks{Krantz's research is supported in part
by NSF Grant DMS-9022140 during residence at MSRI.
Li is partially supported by NSF Grant DMS--9500758.}

\date{August 2, 1995}

\begin{abstract}
Elliptic estimates in Hardy
classes are proved on domains with minimally smooth boundary.
The methodology is different from the original methods of
Chang/Krantz/Stein.
\end{abstract}

\maketitle

\section{Introduction} 
Let $D$ be a bounded Lipschitz domain in
$\RR^n$. In [CKS], Chang, Krantz and Stein introduced certain
real variable Hardy
spaces. They defined two
Hardy spaces $h^p_r(D)$ and $h^p_z(D)$,
$0 < p \leq 1$.  We say that a distribution $f$
lies in $h^p_r(D)$
if it has an extension $E(f)\in h^p(\RR^n)$; and  we say that a
distribution $g$ lies in $h_z^p(D)$ if there is an extension
 $E_z(g) \in h^p(\RR^n)$ so that $E_z(f)=0$ on 
$\RR^n \setminus \overline {D}$.  Here $h^p(\RR^n)$
is a local version,
due to Goldberg [G], of the classical Hardy spaces. 

Since $D$ is bounded,
it is known from Miyachi [M] that
the Hardy space $h_r^p(D)$ can be identified with 
 the subspace of 
distributions $f\in {\cal D}(D)$ such that the radial maximal function
$f^+(x)=f_{\phi, D}^+(x)\in L^p(D)$.  Here
$$
f^+(x)=\sup\Big\{ \Big|\int_D \phi_t(x-y) f(y) dy\Big|: 
0<t< \delta(x) \Big\} ,
$$
where $\delta(x)= \hbox{dis}(x,D^c)$ with
 $D^c=\RR^n\setminus D$ and $\phi$ is a fixed function
such that $\phi\in C^{\infty}_0(B^n)$ ($B^n$ the unit ball in $\RR^n$),
$\phi\ge 0$, and $\int_{\RR^n} \phi(x) dx=1$.  We set 
$\phi_t(x)= t^{-n} \phi(x/t)$.

Let $G(f)$ be the 
solution of  the Dirichlet boundary problem:
$$
\Delta u=f\quad \hbox{in } \ D,\quad \hbox{and } \quad u \biggr |_{\d
D}=0; \leqno(1.1)
$$ 
also let $N(f)$ be the solution of the Neumann boundary value problem:
$$
\Delta u=f\quad \hbox{in } \ D,\quad \hbox{and } 
\quad{\d u \over \d \nu} \biggr |_{\d D}=0 , \leqno(1.2)
$$ 
where $\nu(x)$  denotes the outward unit normal to the boundary $\d D$ at $x$,
 and ${\d \over \d \nu}=\nabla \cdot \nu$.


In [CKS], Chang, Krantz and Stein extended a classical theorem in $\RR^n$
to a smoothly bounded domain in $\RR^n$ and they proved the following theorem.

\begin{theorem} \sl Let $D$ be a bounded domain in $\RR^n$ with
$C^{\infty}$ boundary.  Let $0<p\le 1$. Then 
$$
\left\|{\d^2 G(f)\over \d x_i \d x_j}\right\|_{h^p_r(D)}+
\left\|{\d^2 N(f)\over \d x_i \d x_j}\right\|_{h^p_r(D)}
\le C_p(D) \|f\|_{h_z^p(D)}.
$$
\end{theorem} 

The question of what is the minimum smoothness condition on $\d D$ so
that the above theorem remains true when $D$ is a bounded domain in
$\RR^n$ is still open.  Certainly one may ask: Does the smoothness of
$\d D$ in Theorem 1.1 depends on $p$?  In [CKS], it was conjectured
that Theorem 1.1 remains true if $D$ has $C^{1/p}$ boundary.

The Dirichlet problem for the inhomogeneous Laplace equation (1.1)
has been studied by many authors (see, for example, [AND], [GT], [FKP],
[LiM], and [Ken2]).  Fix $0\le \alpha <\infty$ and $1<p<\infty$.
When $D$ is a smoothly bounded ($C^{\infty}$) domain in $\RR^n$,
Calder\`on-Zygmund theory shows (see [AND]) that if 
$f\in W^p_{\alpha}(D)$ (the Sobolev space),
then there exists a unique $u$ which 
solves (1.1) with
$$
\|u\|_{W^p_{2+\alpha}(D)}\le C\|f\|_{W^p_{\alpha}(D)},
\quad  1<p<\infty.  \leqno(1.3)
$$

In [Dah], B. Dahlberg constructed  a bounded Lipschitz domain 
$D\subset \RR ^2$, and an $f\in C^{\infty}(\Dbar)$ 
so that the regularity (1.3) fails when $\alpha =0$.  In [JeK],
D. Jerison and C. Kenig constructed  a bounded domain in $\RR^2$
with $C^1$ boundary and a function $f\in C^{\infty}(\Dbar)$
so that (1.3) fails for $p=1$ with $\alpha=0$.
Since $C^{\infty}(\Dbar)\subset h^1_z(D)$, we see that the 
aforementioned conjecture is not true when $p=1$. 

>From the definition of $p$-atom (see [CKS]), one can
see that the necessary order of cancellation in an atom
depends not only on $p$ but also
 on $n$. Since the dual of $H^p$ is $\dot{\Lambda}_{n(1/p-1)}$ (the 
homogeneous Zygmund
class), it seems that a reasonable necessary requirement on the smoothness
of $\d D$ so that Theorem 1.1 remains true
 will be that the boundary is  $C^{n(1/p-1)}$ when $p<n/(n+1)$.

The primary purpose of the present paper is to prove that this
last suggested necessary condition is also sufficient. 
We will work
on domains with at least $C^2$ boundary.
As we shall see, our results will be o interest for $p$ small:
$0 < p \leq n/(n+2)$. We shall leave the case when $D$ has $C^q$ boundary
with $1<q\le 2$ for a future paper.

For $0<p<\infty$ and any $\epsilon>0$, we let
$$
\alpha(p,\epsilon)=\max\{2+\epsilon, n(1/p - 1) + \epsilon\}.\leqno(1.4)
$$
The first theorem we propose to prove is:

\begin{theorem} \sl Let $0<p<\infty$ and 
let $D$ be a bounded  domain in $\RR^n$ with $C^{\alpha(p,\epsilon)}$ boundary.
 Then 
$$
\left\|{\d^2 G(f)\over \d x_i \d x_j}\right\|_{h^p_r(D)}+
\left\|{\d^2 N(f)\over \d x_i \d x_j}\right\|_{h^p_r(D)}
\le C_{p,n,\epsilon}(D) \|f\|_{h_z^p(D)}
$$
for any $\epsilon>0$.
\end{theorem} 

\noindent{\bf Note 1:} From the proof of Theorem 1.2, 
it may be seen that we need such an
$\epsilon>0$ only when $n(1/p-1)$ is an integer.
\medskip

We shall show by example that Theorem 1.2 fails for some domain
$D$ having only $C^{n(1/p-1)-\epsilon}$ boundary.

The proof of Theorem 1.1 given by Chang, Krantz and Stein is based
on mapping $D$ to the model domain $\RR^{n}_+$. Our approach for
 proving Theorem 1.2
will be based instead on the machinery connected with the
fundamental properties of the Green's function
for $-\Delta$ in $D$. The properties of the Green's function 
have, historically, played a 
crucial role in
solving the Laplace equation. 
M. Gr\"uter and K.-O. Widman [GW] as well as
 E. Fabes and W. Stroock [FaS] studied the Green's function in a Lipschitz domain
in $\RR^n$. They  gave the basic estimates for the 
Green's function and its first derivative. 
In order to prove Theorem 1.2, we need the asymptotic behavior
of $G$ and its {\em higher derivatives} near the boundary. The
secondary purpose of the present paper is to estimate the
Green's function and its derivatives of all orders (see Theorem 2.2.)
 We believe that the properties of the
Green's function that are derived
in this paper will be useful in  other contexts as well.

\medskip
\noindent {\bf Note 2:}  From our proof of Theorem 1.2, we have that 
$$
{\d^2 G(f)\over \d x_i \d x_j} 
=\Gamma_2(f)+U_2(f),
$$ where
$\Gamma_2$ is bounded from $h^p_z(D)$ to $h_r^p(D)$ for $D$ with
Lipschitz 
boundary.
The operator $U_2: h^p_z(D) \to h^p(D)$ is 
bounded, where $h^p(D)$ is the subspace
of $h_r^p(D)$ consisting of all harmonic functions on $D$.

>From our estimates on the higher derivatives of the Green's function, we
are able to consider the elliptic boundary value problems (1.1) and (1.2)
with $f$ lying in the Hardy-Sobolev space. For this purpose, let us 
now introduce
the definition of Hardy-Sobolev space. 

Let $k$ be a non-negative integer and $0<p <\infty$.
We let $h^{k,p}_z(D)$ denote
the space of all measurable functions $f$ with the weak derivative
$\nabla ^{\ell} f\in h^p_z(D)$
for all $0\le \ell \le k$; here $\nabla^k f$ denotes all $k^{\rm th}$ 
derivatives
of $f$. 
We say that  a measurable function $f$ on $D$ belongs to
$h_r^{k,p}(D)$ if $ \nabla ^{\ell} f \in h^p_r(D)$ for all $0\le \ell\le k$.
It is obvious that $h_z^{k,p}(D)\subset h_r^{k,p}(D)$. We shall prove
the following:

\begin{theorem} \sl Let $0<p \le 1$ and let $k$ be a non-negative integer 
with $k\le n(1/p-1)$.
Let $D$ be a bounded  domain in $\RR^n$ with $C^{3+[n(1/p-1)]}$ boundary. Then 
$$
\left\|{\d^{2} G(f)\over \d x_i \d x_j }\right\|_{h^{k,p}_r(D)}+
\left\|{\d^2 N(f)\over \d x_i \d x_j }\right\|_{h^{k, p}_r(D)}
\le C_{p,k,n}(D) \|f\|_{h_z^{k,p}(D)}.
$$
\end{theorem} 

The paper is organized as follows. In Section 2, we estimate the behavior
of the higher derivatives of the Green's function near the boundary. The
results in this section are the key to the rest
of the paper. In Section 3, we use the
atomic decomposition theorem in [CKS] to reduce the proof of Theorem 1.2
to the study of a single atom.  In Section 4, we prove 
Theorems 1.2 and 1.3 for the Dirichlet boundary value
problem. In Section 5, Theorems 1.2  and 1.3 for the Neumann boundary 
problem are proved and the
`Green's function' for the Neumann boundary value problem is studied.
In Section 6, we shall construct examples that show that 
Theorem 1.2 is reasonably
sharp.

\medskip

The authors wish to thank Jiaping Wang for a useful conversation
on the Green's function.
\medskip

\section{Estimate derivatives of Green function}  
It is well-known that
 the fundamental solution for the Laplacian $\Delta$ in $\RR^n$ is
$$
\Gamma(x-y)=\Gamma(|x-y|)=
\cases{  -(n-2)^{-1} \omega_n^{-1}  |x-y|^{2-n},  & if $n>2$\cr
{1\over 2\pi } \log |x-y|, & if $n=2$\cr}
$$
where $\omega_n$ denotes the surface area of the unit sphere in
$\RR^n$ (see [KR]).

Let $D$ be a bounded $C^1$ 
domain in $\RR^n$.  
For each $y\in D$, we let $U(\cdot, y)$
be the solution of the Dirichlet problem:
$$
\Delta U(\cdot, y)=0 \hbox{ in } D,\quad U(x,y)=\Gamma(y-x),
 \hbox{ for  } x\in 
\d D. \leqno(2.1)
$$
Then the Green's function for $\Delta$ on $D$ is
$$
G(x, y) =\Gamma(x-y)-U(x,y).   \leqno(2.2)
$$
The Dirichlet problem (1.1) 
has a unique solution
$$
u(x)=G(f)(x)=\int_D G(x, y) f(y) \, dy.   \leqno(2.3)
$$

Now let $V(\cdot, y)$  be the solution of the following Neumann problem:
$$
\Delta V =\hbox{Vol}(D)^{-1}\ \  \hbox{ in } \ D,\quad 
{\d \over \d \nu}V={\d \over \d \nu} \Gamma (y-\cdot)
 \ \ \hbox{ on } \   
\d D. \leqno(2.4)
$$
We let
$$
N(x,y)=\Gamma (x-y)-V(x, y)  \leqno(2.5)
$$
It is easy to show that the Neumann problem (1.2) has a unique solution,
up to an additive constant, given by the formula
$$
u(x)=N(f)(x)=\int_D N(x,y) f(y) d y.     \leqno(2.6)
$$

The main purpose of this section is to study 
the basic properties of
$U(x,y)$. We shall derive information about the
asymptotic behavior of $U(x,y)$ 
when $x,$ $y$ are near the boundary
$\d D$. For convenience, we will always assume that $n>2$ (the case $n=2$
is similar, but the details of the formula are different.) 
Similar results for $V(x,y)$
will be obtained in Section 5.

We will need the following proposition that is due to  M. Gr\"uter 
and K.-O. Widman [GW].

\begin{proposition}\sl
 Let $D$ be a bounded  domain in $\RR^n$ satisfying the uniform
exterior sphere condition (i.e., each boundary point of $D$
has an exterior osculating sphere of uniform size). 
Then  the Poisson kernel
$$
P(x,y)={\d G(x, y)\over \d \nu(y)}
$$
satisfies the estimate:
$$
0\le P(x,y)\le {C\delta (x)\over |x-y|^n}\leqno(2.7)
$$
for all $x\in D$ and $y\in \d D$.  {\sl [Note that, on a $C^2$ domain,
these last two expressions are known to be comparable---see [KR]]}. 
\end{proposition}

The main result of this section is to prove the following theorem about
$U$. 

\begin{theorem}\sl Let  $k_0\ge 2$ be an
integer. Let
 $D$ be a bounded domain with
$C^{k_0,\epsilon} $ boundary. Then for any
multi-indices $\alpha$ and 
$\beta$ with $|\alpha|=k\le k_0$ and $|\beta|=\ell\le k_0$,
 we have
$$
\left|{\d^{k+\ell} U(x, y)\over \d x^{\beta}  \d y^{\alpha}}\right|
\le {C_{k_0,\epsilon}\over (|x-y|+\delta(y))^{n +k+ \ell -2}}
\leqno(2.8)
$$
for all $x,\ y \in D$ and any $\epsilon>0$. 
\end{theorem}

\proof  If we write $f(t)=- t^{1-n/2}$ for all $t\ge 0$, then
$\Gamma(x-y)=f(|x-y|^2)$. Now 
$$
{\d^k U(x,y) \over \d y^{\alpha}}=\int_{\d D} P(x, z) {\d^k \Gamma(z-y)
\over \d y^{\alpha}} \, d\sigma(z) ,
\leqno(2.9)
$$
where $P(x,z)$ is the Poisson kernel satisfying (2.7) for all $x\in D$ and
$z\in \d D$ (since $D$ has at least $C^2$ boundary)  Thus
$$
{\d^k U(x, y)\over \d y^{\alpha}}=
{\d^k \Gamma(x-y)\over \d y^{\alpha}}\ , 
     \qquad x\in \d D,\  y\in D.  \leqno(2.10)
$$

For any fixed
 $x_0,\,  y\in D$ , we let $R=|x_0-y|/2$. If $R<16 \delta(y)$, then
(2.8) holds with $(x,y)=(x_0,y)$ by the maximum principle and Proposition 2.1.
We now assume that $R\geq 16 \delta(y)$, and we let 
$\epsilon_0=\delta(x_0)^{n+\ell+k} \delta(y)^{n+2k}$. 
For any multi-index $\alpha$ with $|\alpha|=k$, we consider the function 
$$
f_{\alpha, y}(|x-y|^2+\epsilon_0)
={\d^k f(|x-z|^2+ \epsilon_0) \over \d z^{\alpha}
}\biggr |_{z=y} .
$$ 
Since $U(\cdot,y)$ is harmonic, we have
\begin{eqnarray*}
\lefteqn { {\d ^{\ell+k} U(x, y)\over \d x^{\beta} \d y^{\alpha}}-
{\d^{|\beta|} f_{\alpha, y}(|x-y|^2+ \epsilon_0)\over \d x^{\beta}}}\\
&=& \int_{\d D}{\d^{\ell}  P(x, z)\over \d x^{\beta}}
 \Big({\d ^k \Gamma(y-z)\over \d y^{\alpha}}-
f_{\alpha, y}(|y-z|^2+\epsilon_0 )\Big) \, d\sigma(z)\\
&& -{\d^{\ell}\over \d x^{\beta}} \int_D  G(x, z)
 \Delta_z f_{\alpha, y}(|y-z|^2+\epsilon_0) \, dz\\ 
& \equiv & J_1(x,y) + J_2(x,y)
\end{eqnarray*}
Notice that the term involving integration over the interior comes from Green's
identity---since
the integrand has a singularity.

Since $ \d D$ is $C^{k_0,\epsilon}$ we have (from the maximum principle) that
$$
\biggl |{\d^{|\gamma|} U(x, z)\over \d z^{\gamma}}\biggr |
=\biggl | {\d^{|\gamma|} U(z, x)\over \d z^{\gamma}}\biggr | 
\le 
 C_{\gamma,\epsilon}\|\Gamma(x-\cdot)\|_{C^{|\gamma|} (\d D)}
\le C_{n,\gamma} \delta (x)^{-|\gamma|-n+2}
$$
for all $|\gamma|\le k_0$. Moreover, since again $U(x,\cdot)$ is
harmonic, we have
$$
\biggl |{\d^{|\gamma|+|\beta|} G(x, z)\over \d x^{\beta} \d z^{\gamma}
}\biggr |
\le C_{n,\gamma, \epsilon,\beta}
\Big(|x-z|^{-n-|\gamma|-|\beta|+2}+
 \delta (x)^{-|\gamma|-|\beta|-n+2}\Big)  \leqno(2.11)
$$
for all $|\gamma|\le k_0$ and any $\beta$.

Since $P(\cdot, z)$ (for $z\in \d D$) is harmonic,
we have 
$$
\biggl |{\d^{\ell} \over \d x^{\beta}} P(x,z) \biggr |
\le C_{\ell} \delta (x)^{-n-\ell+1}
$$
for all $z\in \d D$. Thus
\begin{eqnarray*}
|J_1(x_0,y)|&\le&
 C_{k,\ell}\int_{\d D} \delta(x_0)^{-\ell-n} {\delta(x_0)^{\ell+k+n}
\delta(y)^{n+2k}  \over
(|z-y|^2 +\epsilon_0 )^{(n+k-1)/2}} \, d\sigma(z) \\
&\le& C_{k,\ell} \int_{\d D} {\delta(y)^{2k+n}\over \delta(y)^{n+k-1} }
\, d\sigma(z) \\
&\le& C_{k,\ell}.
\end{eqnarray*}

Also, for $R =|x_0-y|/2 > 0$ fixed,
\begin{eqnarray*}
|J_2(x_0,y)|&\le & \Biggl | {\d ^{\ell} \over \d x^{\beta}}
\int_{D\setminus B(y, R)} G(x,z)
 \Delta_z f_{\alpha, y}(|y-z|^2+\epsilon_0) \, d z\biggr |_{x=x_0}\Biggr |\\
&& +\biggl |\int_{D\cap B(y, R)}{\d^{\ell}  G(x_0, z) \over \d x^{\beta}}
\Delta_z f_{\alpha, y}(|y-z|^2+\epsilon_0) \, d z\biggr |\\
&\equiv& J_{21}(x_0,y)+ J_{22}(x_0,y).
\end{eqnarray*}
Now
\begin{eqnarray*}
J_{21}(x_0,y)
&\le& 
 \Biggl | {\d ^{\ell}\over \d x^{\beta}}\int_{D\setminus B(y, R)}
\Gamma(x-z) \Delta_z f_{\alpha, y}(|y-z|^2+\epsilon_0)
\, d z \biggr |_{x=x_0}\Biggr |\\
& &+ \int_{D\setminus B(y, R)}\epsilon_0
\delta(x_0)^{-n-\ell+2} (|z-x_0|+R)^{-n-k} \, d z\\
&\le&J_{211}+ C_{k,\ell} |x_0-y|^{-n-k-\ell +2} ,
\end{eqnarray*}
where $J_{211}$ is defined by the last inequality.
We see that 
\begin{eqnarray*}
J_{211}(x_0,y)
&=&
 \Biggl | {\d ^{\ell}\over \d x^{\beta}}\int_{D\setminus B(y, R)}
\Gamma(x-z) \Delta_z f_{\alpha, y}(|y-z|^2+\epsilon_0)\,
d z\biggr |_{x=x_0}\Biggr |\\
&\le&
 \Biggl | {\d ^{\ell}\over \d x^{\beta}}\int_{\d(D\setminus B(y, R))}
\Gamma(x-z) D_{\nu(z)} f_{\alpha, y}(|y-z|^2+\epsilon_0)
\, d\sigma (z)\biggr |_{x=x_0}\Biggr |\\
&&+
 \biggl | {\d ^{\ell}\over \d x^{\beta}}\int_{\d(D\setminus B(y, R))}
D_{\nu(z)}\Gamma(x-z) f_{\alpha, y}(|y-z|^2+\epsilon_0)
\, d\sigma(z)|_{x=x_0}\biggr |\\
&&+ \Biggl | {\d^{\ell} f_{\alpha, y}(x,y)\over \d x^{\beta}} 
   \biggr |_{x=x_0} \Biggr | \\
& \leq & C_{k,\ell} |x_0 - y|^{-n-k-\ell+2}
\end{eqnarray*}
by standard arguments (since $k\le k_0$,
$\d D$ is $C^{k_0,\epsilon}$ with $k_0\ge 2, \ \epsilon>0$
and by the fact that $|\nu(z) \cdot (x-z)|\le C(|x-z|^2+\delta(x)$
for all $z\in \d D$ and $x\in D$).  [Note that the last 
term in the penultimate line comes from Green's theorem.]

Now we consider $J_{22}(x_0,y)$. Let $\Omega= D\cap B(y, R)$. 
After applying the divergence theorem many times, we have
\begin{eqnarray*}
\lefteqn{J_{22}(x_0,y)}\\
&=&
\Big|\int_{\Omega }{\d^{\ell} G(x_0, z)\over \d x^{\beta}}
 \Delta_z f_{\alpha,y}(|y-z|^2+\epsilon_0) d z\Big|\\
&=& \Big|\int_{\Omega}{\d^{\ell}  G(x_0,z)\over \d x^{\beta}} 
{\d^k \over \d z^{\alpha}}\Delta_z f(|z-y|^2+\epsilon_0) d z\Big|\\
&\le& \Big|\int_{\Omega} {\d^{\ell+k} G(x_0,z)\over \d x^{\beta}
\d z^{\alpha}} 
\Delta_z f(|z-y|^2+\epsilon_0) d z\Big| \\
&&+ \Big|\int_{\d \Omega}\sum_{|\gamma|\le k}|(\nu(z)^{\beta}| 
\biggl |{\d^{\ell+|\gamma|}G(x_0,z)\over \d x^{\beta}\d
z^{\gamma}}\biggr | 
\biggl |{\d^{k-|\gamma|} \over \d z^{\alpha-\gamma} } \Delta_z 
f(|z-y|^2+\epsilon_0)\biggr | d \sigma(z)\\
& \equiv & J_{221}(x_0,y)+ J_{222}(x_0,y).
\end{eqnarray*}
Now, since $\Delta_z f(|y-z|^2)=0$ for all $z\ne y$, we have
$$
\Big|\Delta_z f(|y-z|^2+\epsilon_0)\Big|={ n(n-2)\,  \epsilon_0
\over (|z-y|^2+\epsilon_0)^{n/2}}
\le C_{n,k}{\delta(x_0)^{n+k+\ell}\delta (y)^{2k+n} \over
|z-y|^n}\leqno(2.12) .
$$
Moreover, we have 
$$
\Big|D^{\gamma}_z \Delta_z f(|y-z|^2+\epsilon_0)\Big|
\le C_{n,\gamma}{\delta(x_0)^{k+\ell +n} \delta(y)^{n+2k} 
 \over (|z-y|+\epsilon_0)^{n+|\gamma|}}
$$
Applying (2.12), we find that 
\begin{eqnarray*}
\lefteqn{J_{222}(x_0,y)}\\
&\le& C_{\ell,k}\int_{\d \Omega}
 \delta(x_0)^{-\ell-k-n+2}{\delta(x_0)^{\ell+n+k}
 \delta(y)^{2k+n}\over
\delta(y)^{n+k}} d\sigma(z)\\
&\le& C_{n,k,\ell}.
\end{eqnarray*}

Finally, we estimate $J_{221}(x_0,y)$. Observe that
\begin{eqnarray*}
\lefteqn{J_{221}(x_0,y) }\\
&=& \Big|\int_{\Omega} {\d^{k+\ell} G(x_0,z)\over \d x^{\beta}
 \d z^{\alpha } }
\Delta_z f(|y-z|^2+\epsilon_0) d z \Big|\\
&\le &  \int_{\Omega} C_{k,n} \delta(x_0)^{-k-\ell-n+2}
{\delta(x_0)^{k+n+\ell} \delta(y)^{n+2k}\over (|y-z|+\epsilon_0)^n} d z\\
&\le& C_{n,k} \, \epsilon_0^{1\over n+k+\ell}\, \biggl (\log{1\over
\epsilon_0}
\biggr )\\
&\le& C_{n,k}.
\end{eqnarray*}
Combining the above estimates, we have
$$
\biggl | {\d ^{k+\ell} U(x_0,y) \over \d x^{\beta} \d y^{\alpha}}
-{\d^{\ell} f_{\alpha}(|x_0-y|^2+\epsilon_0)\over \d x^{\beta}}\biggr |
\le C_{n,k} (\delta(y)+|x-y|)^{-n-k-\ell+2}
$$

Using this estimate and the definition of $f_{\alpha,y}(x_0)$,
we conclude that the proof of Theorem 2.2 is complete.\epf

\begin{corollary}\sl Let $\lambda> 2$ be a non-integer
and   let  $D$ be a bounded domain with
$\lambda $ boundary. Then  for all $0\le k\le [\lambda]$ 
and  any cube $Q$ in $\RR^n$ with $2Q\subset D$, we have
$$
\left|{\d^2 U(x, y)\over \d x_i \d x_j}-\sum_{|\alpha|\le k}
{\d^{\alpha}\over \d y^{\alpha}}{\d^2 U(x, x_0)
\over \d x_i \d x_j} (y-x_0)^{\alpha}\right|
\le C_{p,n,\lambda} {\delta^{k(\lambda)}\over |x-x_0|^{n+k(\lambda)}}
\leqno(2.13)
$$
for all $x\in D\setminus 2Q$ and $y, x_0\in Q$ and $|Q|=\delta^n$. Where $k(\lambda)=k+1$ if $k<[\lambda]$, and $k(\lambda)=\lambda$ if $k=[\lambda]$ 
\end{corollary}
\proof This follows directly from Theorem 2.2 
(with a suitable modification for
fractional derivatives) and from Taylor's theorem.

\section{Reduction of Theorem 1.2}

Let $f$ be a measurable function on $D$ and let $\lambda$ and $q$ be
positive numbers.
 Let us recall the following definition for a $p$-atom.
 Let $a$ be a bounded
function in $\RR^n$; we say that $a$ is an $h^p(\RR^n)$ atom
if $a$ is supported in a cube $Q$ with  $(\int_{\RR^n} |a(x)|^2
dx )^{1/2} \le |Q|^{1/2-1/p}$ and either {\bf (i)}
$|Q|>1$ or {\bf (ii)} $|Q| \leq 1$ 
and, for each $\alpha=(\alpha_1,\cdots, \alpha_n)$ with
 $|\alpha|\le [n(1/p-1)]$ , we have
$$
\int_Q x^{\alpha} a(x) d x=0.
$$
where all $\alpha_i\ge 0$ and $|\alpha |=
\sum_{i=1}^n \alpha_i$.

In order to prove  Theorem 1.2,
we recall the following theorem that is formulated from 
the statement of Theorem 3.2 in [CKS] and its proof.

\begin{theorem}\sl Let $D\subset \RR^n$ be a bounded Lipschitz domain and
let $0<p\le 1$. Then 
$f\in h^p_z(D)$ if and only if 
$f$ has an atomic decomposition, with $h^p(\RR^n)$ atoms whose supports
lie in $D$, i.e.,
$$
f=\sum_{Q\subset D} \lambda_Q a_Q
$$
with $2Q\cap \d D=\emptyset$ if the diameter of $Q$ is small (i.e.\ 
$\leq 1$), and
$$
\sum_{Q\subset D} |\lambda_Q|^p<\infty.
$$
\end{theorem}

By Theorem 3.1, it is easy to verify that 
Theorem 1.2 holds  if and only if
it  holds  for $f$ an $h^p$-atom. 
More precisely, the proof of Theorem 1.2 can be reduced to proving
the following theorem:

\begin{theorem} \sl Let $0<p\le 1$.
 Let $D$ be a bounded  domain in $\RR^n$ with 
$C^{\alpha(p,\epsilon)}$ boundary.
 Then 
$$
\left\|{\d^2 G(a)\over \d x_i \d x_j}\right\|_{h^p_r(D)}+
\left\|{\d^2 N(a)\over \d x_i \d x_j}\right\|_{h^p_r(D)}
\le C_{p,n,\epsilon}(D).
$$
holds for all $p$-atoms $a$ with support $Q$ satisfying $2Q\subset D$.
\end{theorem}

We shall separate the proof of Theorem 3.2 and Theorem 1.3
 into two cases, which are
given in Sections 4 and 5 respectively.

\section{The Dirichlet Problem}

In this section, we shall prove Theorem 3.2 with the Dirichlet
 boundary condition. In other words, we shall prove Theorem 4.1.

\begin{theorem} \sl Let $0<p\le 1$ and let
 $D$ be a bounded  domain in $\RR^n$ with $C^{\alpha(p,\epsilon)}$ boundary.
 Then 
$$
\left\|{\d^2 G(a)\over \d x_i \d x_j}\right\|_{h^p_r(D)}
\le C_{p, n,\epsilon} (D).  \leqno(4.1)
$$
holds for all $p$-atoms with support $Q$ satisfying $3Q\subset D$ and
for any $\epsilon>0$.
\end{theorem} 

Since
$$
G(a)(x)=\int_{D} \bigl (\Gamma(x-y)-U(y,x) \bigr )a(y) \, dy ,
$$
it is clear that 
$$
\biggl \|{\d ^2\over \d x_i \d x_j}\int_D \Gamma(y-x) a(y) 
\, d y \biggr \|_{h_r^p(D)}\le C_{p,n}
$$
from the results in [ST1] for Hardy spaces in $\RR^n$. 
Let $a$ be an $h^p$-atom.  Without loss of generality, we may assume that
$0\in D$ and 
$$
\supp(a) \subset \{ y \in D: |y| < \delta\}=B(0,\delta)\subset D.
$$
If the atom $a$ is supported in a cube $Q\subset D$ with diameter 
greater than
$1$, then $a\in L^2(\RR^n)$. Then  we have that (4.1) 
holds by using a modified
version of the
following result of Jerison and Kenig in [JeK] and Gilbarg and Trudinger [GT].

\begin{theorem} \sl Let $D$ be a bounded  domain
in $\RR^n$ with $C^2 $ boundary. If $u$ is a solution
of (1.1), then
$$
\|u\|_{W^p_{2+\alpha}(D)}\le C_{p,\alpha}\|f\|_{W^p_{\alpha}(D)}.
$$ 
for all $1<p<\infty$ and $-1<\alpha\le 0$.
\end{theorem}

\noindent Note, in passing, that there is certainly an analogous
version
of Theorem 4.2 for the Neumann problem.

Now we assume that $a$ is a classical atom with support in a cube $Q$ with 
small diameter with $2Q\cap \d D= \emptyset$. Let the center of $Q$ be
$x_0=0$. For each $x\in D$, we let
$$
H(a)(x)={\d^2 \over \d x_i \d x_j}\int_D U(x, y) a(y) dy.
$$
Since $U(x,y)$ is harmonic in both
 $x$ and $y$, we have that $H(x)$ is harmonic in $x$. 
 We choose $\phi\in C^{\infty}_0(B(0,1))$ to be a non-negative
radial function such that 
 $\int_{\RR^n} \phi(x) dx=1$. Let $x\in D$ and 
 $d(x)=\hbox{dis}(x, D^c)$. Since $H(a)$ is harmonic and $\phi$
is radial, the mean value property of harmonic function
shows that
$$
H(a)^+(x)=\sup_{0<t < d(x)} \left| \int_{D} 
H(a)(y) \phi_t(x-y) dy\right|=|H(a)(x)|,
$$
where $\phi_t=t^{-n} \phi( t^{-1} \, x)$.  
Now 
$$
\int_{D}H(a)^+(x)^p dx
=\int_{B(0,4\delta)}|H(a)(x)|^p dx
+ \int_{D \setminus B(0, 4\delta)}|H(a)(x)|^p dx
 \equiv I_1(a)+I_2(a).
$$
By (1.3), we have
\begin{eqnarray*}
I_1(a)&\le & \left(\int_{B(0,4\delta)} |H(a)(x)| ^2 d x\right )^{p/2} 
\left(\int_{B(0,4\delta)} dx\right) ^{(2-p)/2}\\
&\le & C^p \left(\int_{B(0,4\delta)} |a(x)| ^2 d x\right)^{p/2} 
|B(0,4\delta)|^{(2-p)/2}\\
&\le & C^p |B(0,4\delta)|^{(-2/p+1)p/2} 
|B(0,4\delta)|^{(2-p)/2}\\
&=& C^p.
\end{eqnarray*}

Next we  estimate $I_2(a)$.  By Corollary 2.3, we have
$$
\left|{\d^2 U(x, y)\over \d x_i \d x_j}-\sum_{|\alpha|\le n_p}
{\d^{\alpha}\over \d y^{\alpha}}{\d^2 U(x, x_0)
\over \d x_i \d x_j} (y-x_0)^{\alpha}\right|
\le C_p {\delta^{n(1/p-1)+\epsilon/2}\over |x-x_0|^{n+n(1/p-1)+\epsilon/2}}
$$
for all $x\in D\setminus 4B(x_0, \delta)$ and $ y\in 
B(x_0, 2\delta)\subset D$,
where $n_p=[n(1/p-1)]$.
Thus, since the center of $Q$ is $x_0=0$, for each 
$x\in D\setminus B(0, 4\delta)$ we have
\begin{eqnarray*}
|H(a)(x)| &=& \Big|\int_Q {\d^2 U(x,y)\over \d x_i \d x_j} a(y) \, dy\Big|\\
&=&  \biggl |\int_Q
\Big ({\d^2 U(x, y)\over \d x_i \d x_j}-\sum_{|\alpha|\le n_p}
{\d^{\alpha}\over \d y^{\alpha}}{\d^2 U(x, 0)
\over \d x_i \d x_j} y^{\alpha}\Big ) a(y) \, dy\biggr | \\
&\le&  C_{p,\epsilon}\int_Q  {\delta^{n(1/p-1)+\epsilon/2}\over
|x-x_0|^{n+n(1/p-1)+\epsilon/2}} |a(y)| \, dy \\
&\le& C_{p,\epsilon} |Q|^{1-1/p}  { \delta^{n(1/p-1)+\epsilon/2} \over |x|^{n+n(1/p-1)+\epsilon/2} }\\
&= & C_{p ,\epsilon} { \delta^{\epsilon/2}\over |x|^{n/p+ \epsilon/2} }.
\end{eqnarray*}
Therefore, if $D \subseteq B(0,d_0),$ then
\begin{eqnarray*}
\int_{D\setminus B(0, 4\delta)}|H(a)(x)|^p dx 
&\le& C_{p,\epsilon}^p\int_{D\setminus B(0, 4\delta)}  
{ \delta^{p\epsilon/2 }\over |x|^{n+p\epsilon/2}} dx\\
&\le& C_p^p \int_{4\delta}^{d_0} C_n  
{ \delta^{p\epsilon/2}\over t^{1+ p\epsilon/2 }} dt\\
&=& C_p^p {2\over p\epsilon}  { \delta^{p\epsilon/2}\over 
(4 \delta)^{p\epsilon/2}}\\
&= & C(p,n,d_0).
\end{eqnarray*}
The proof of Theorem 4.1 is thus complete.\epf

\begin{theorem} \sl Let $0<p <\infty$ and let $k\le n(1/p-1)$ be
a non-negative integer.  Let
 $D$ be a bounded  domain in $\RR^n$ with $C^{2+n(1/p-1)+\epsilon}$ boundary  
for some $\epsilon > 0$.  Then 
$$
\left\| {\d ^2 G(f) \over \d x_i \d x_j} \right\|_{h^{k,p}_r(D)}
\le C_{k,p,\epsilon} \|f\|_{h^{k,p}_z(D)}.  \leqno(4.2)
$$
for any $\epsilon>0$.
\end{theorem}

\proof Observe that
$$
G(f)(x)=\int_D \Gamma(x-y) f(y) dy- \int_D U(x,y) f(y) dy
=\Gamma (f)(x)+ U(f)(x)
$$
Since $f$ has compact support, it is obvious that 
$$
\left\| \nabla^2 \Gamma (f) \right\|_{h^p_r(D)}
\le C_{k,p,\epsilon} \|f\|_{h^{k,p}_z(D)}.  \leqno(4.3)
$$

Next we prove that
$$
\left\| \nabla^2 U(f) \right\|_{h^p_r(D)}
\le C_{k,p,\epsilon} \|f\|_{h^{k,p}_z(D)}.  \leqno(4.4)
$$
Now $\nabla^{\ell} U(a)(x)$ is harmonic for any 
non-negative integer $\ell$, so it suffices to show 
that $\|\nabla^{\ell} U(a)\|_{L^p(D)}\le C_p $ for all
$0\le \ell \le k+2$. We shall prove the case $\ell=k+2$;
the other cases are even easier. 
In order to do this, we need the following Sobolev 
embedding theorem which is a special case of Theorem 2 in [BB]
and can also be deduced from Theorem 2 in [HPW]:
$$
h^{k,p}_z(D)\subset h_z^{pn/(n-pk)}(D)\leqno(4.5)
$$
and the embedding is continuous.
\medskip
By Theorem 3.1 and the fact that $\nabla^{\ell} U(f)$ is harmonic, it 
suffices to prove that 
$$
\left\|{\d^{k+2} U(a) \over \d x^{\alpha}}\right \|_{L^p(D)} \le C_{k,p,\epsilon}
$$
for any $pn/(n-pk)$--atom with support $Q$ and $3Q\cap \d D=\emptyset$
(since $k\le n(1/p-1)$, hence $pn/(n-pk)\le 1$.)

Notice that
$$
{\d^{k+2} U(a)(x) \over \d x^{\alpha}}
=\int_{\Omega} {\d^{k+2} U(x,y) \over \d x^{\alpha}} a(y) dy .
$$
Let $x_0$ be the center of $Q$. Then
\begin{eqnarray*}
\lefteqn{ \int_{2Q}{\d^{k+2} U(a) \over \d x^{\alpha}}(x)|^p dx}\\
&\le & C_{p,k, n,\epsilon}\int_{2Q} \delta (x)^{-kp-np} |Q|^{-(n-pk)/n +p}  dx\\
&\le & C_{p,k,n,\epsilon}\int_{2Q} \delta(x_0)^{-kp-np}
 \delta(x_0)^{-n+pk+pn} dx\\
&=&C_{p,k,n,\epsilon}.
\end{eqnarray*}

We set 
$$
\ell_p= [n(n-pk)/np-n)]=[n/p-k-n], \quad \eta(k)=1 \hbox{ if } k>0,\ \ \eta(0)=\epsilon.
$$
Then 
$$
\left|{\d^{k+2} U(x, y)\over \d x^{\alpha}}-\sum_{|\alpha|\le \ell_p}
{\d^{\alpha}\over \d y^{\alpha}}{\d^{k+2} U(x, x_0)
\over \d x^{\alpha} } (y-x_0)^{\alpha}\right|
\le C_p {\delta^{(n/p-k)-n+\eta(k)}\over |x-x_0|^{(n/p-k)+\eta(k)+k}}
$$
for all $x\in D\setminus 2Q$ and $ y\in Q$ as above.
Thus, for any
$x\in D\setminus 2Q $ we have
\begin{eqnarray*}
\lefteqn{\Big|{\d^{k+2} U(a)(x) \over \d x^{\alpha}}\Big|}\\
&=& \Big|\int_Q {\d^{k+2} U(x,y)\over \d x^{\alpha} } a(y) dy\Big|\\
&=&  \biggl |\int_Q
\Big ({\d^{k+2} U(x, y)\over \d x^{\alpha} }-\sum_{|\beta|\le \ell_p}
{\d^{|\beta|}\over \d y^{\beta}}{\d^{k+2} U(x, x_0)
\over \d x^{\alpha}} (y-x_0)^{\beta}\Big ) a(y) dy\biggr | \\
&\le&  C_{p,k,n,\epsilon}\int_Q  
{\delta^{(n/p-k)-n+\eta(k) }\over |x-x_0|^{k+(n/p-k)+\eta(k)}} |a(y)| dy \\
&\le& C_{p,k,n,\epsilon} |Q|^{1-(n-pk)/pn}  { \delta^{(n/p-k)-n+\eta(k)} 
\over |x-x_0|^{k +(n/p-k)+\eta(k)} }\\
&= & C_{p,k,n,\epsilon} { \delta^{\eta(k)}\over |x-x_0|^{k+(n/p-k)+\eta(k)} }
\end{eqnarray*}
for any $\epsilon>0$.
Therefore, if $D \subseteq B(0,d_0),$ then
\begin{eqnarray*}
\lefteqn{\int_{D\setminus 2Q }|H(a)(x)|^p dx }\\
&\le& C_{p,k,n,\epsilon}^p\int_{D\setminus 2Q} { \delta^{p\eta(k)}
\over |x-x_0|^{pk+n -pk+p\eta(k)} } dx\\
 &\le& C_{p,k,n,\epsilon,d_0}^p \delta^{p\eta(k) -p\eta(k)}\\
&=& C_{p,k,n,\epsilon}.
\end{eqnarray*}
The  proof of Theorem 4.3 is therefore complete.\epf

\section{The Neumann Problem}

In this section, we shall prove Theorems 2.2 and 1.3
with the Neumann boundary condition.
More precisely, we shall first prove the following theorem.

\begin{theorem}\sl Let $0<p\le 1$ and let $D$ be a bounded domain with
$C^{\alpha(p,\epsilon)}$ boundary. Then
$$
\left\|{\d N(a)\over \d x_i \d x_j}\right\|_{h^p_r(D)}\le C_p(D)
\leqno(5.1)
$$
for all $h^p(\RR^n)$ atoms with support $Q\subset 2Q \subset D$.
\end{theorem}

\proof Let $a$ be an $h_z^p$ atom with support $Q$ and 
$|Q|=\delta^n$. If $\delta\ge 1$, then $a\in L^2(\RR^n)$. By the version of
Theorem 4.2 that holds for the Neumann problem, we have that (5.1) holds.

Now we assume that $a$ is a classical atom with support $Q$ and
$2Q\cap \d D=\emptyset$. In this case, (4.1) follows from the
argument of the proof of Theorem 4.1 and the following two results
(Theorems 5.2 and 5.3.)

\begin{theorem}\sl Let $k_0\ge 2$ be any positive integer and let
 $D$ be a bounded domain with
$C^{k_0,\epsilon} $ boundary. Then, for any 
multi-indices $\alpha$ 
and $\beta$  with $|\alpha|=k\le k_0$ and  $|\beta|=\ell \le k_0$, we have
$$
\left|{\d^{k+\ell} V(x, y)\over \d x^{\beta}  \d y^{\alpha}}\right|
\le {C_{k_0,n, \epsilon}\over (|x-y|+\delta(y))^{n +k+\ell-2}}
\leqno(5.2)
$$
for all $x,\ y \in D$ any $\epsilon>0$. 
\end{theorem}

\proof  We write $f(t)=- t^{1-n/2}$ for all $t\ge 0$. Then
$\Gamma(x-y)=f(|x-y|^2)$. Let $C_D=\hbox{Vol}(D)^{-1}$. Then 
\begin{eqnarray*}
V(x,y)&=&\int_{\d D} P(x,z) V(z,y) d\sigma(z)
+\int_D G(x, z) C_D d z\\
&=& \int_{\d D} \biggl ( {\d \Gamma(x-z)\over \d \nu(z)}-
{\d U(x, z)\over \d \nu(z)} \biggr )
 V(z,y) d\sigma(z)
+\int_D G(x, z) C_D d z\\
&=& \int_{\d D} {\d \Gamma(x-z)\over d\nu(z)}V(z, y) d \sigma(z)\\
&& \quad -\int_{\d D}{\d U(x, z)\over \nu(z)} V(z,y) d\sigma(z)
+\int_D G(x, z) C_D d z\\
&=& I_1(x,y)-I_2(x,y) +I_3(y)
\end{eqnarray*}

It is obvious that $I_3(y)$ satisfies (5.2) by using estimate (2.8).
Now we consider $I_2(x,y)$. Notice that
\begin{eqnarray*}
I_2(x,y)&=& \int_{\d D}{\d U(x, z)\over \d \nu(z)}) V(z,y) d\sigma(z)\\
&=& \int_{\d D} U(x, z) {\d V(z, y) \over \d \nu(z)} d\sigma(z)\\
&&+ \int_D \Delta_z U(x, z) V(z,y) d z-\int_D U(x,z) \Delta_z V(z,y) d z\\
&=& \int_{\d D} U(x, z) {\d \Gamma(z-y) \over \d \nu(z)} d\sigma(z)
+0 -\int_D U(x,z) C_D d z\\
&=& \int_{\d D} \Gamma(x-z) {\d \Gamma(z-y)\over \d \nu(z)} d\sigma(z)
-C_D\int_D U(x,z) d z.
\end{eqnarray*}
Using arguments similar to those in Section 2, 
we see that $I_2(x,y)$ 
satisfies the estimate (5.2).

For any $x_0, y\in D$, if $|x_0-y|< 4\delta(y)$, then it is easy to prove
that (5.2) holds by replacing $(x,y)$ by $(x_0,y)$. Without loss of generality,
we may assume that $|x_0-y|>4\delta(y)$.
 We consider the term $I_1(x,y)$. If we let
$$
f_{x,0}(|z-x|^2+\epsilon_0)=(|x-z|^2+\epsilon_0)^{1-n/2},\quad \epsilon_0=
\delta(x_0)^{n+\ell+k} \delta(y)^{n+2k},
$$ then
\begin{eqnarray*}
I_1(x,y)&=& \int_{\d D}{\d \Gamma(x-z)\over d\nu(z)}V(z, y) d \sigma(z)\\
&=& \int_{\d D} {\d f_{x, 0}(|z-x|^2+\epsilon_0) \over \d \nu(z)}
V(z,y) d\sigma(z)\\
&& + \int_{\d D}\Big({\d \Gamma(x-z)\over \d\nu(z)}-
 {\d f_{x, 0}(|z-x|^2+\epsilon_0)\over \d \nu(z)}\Big)V(z, y) d \sigma(z)\\
&=& I_{11}(x,y) +I_{12}(x,y)
\end{eqnarray*}
We have that
\begin{eqnarray*} 
I_{11}(x,y)&=& \int_{\d D}  f_{x, 0}(|z-x|^2+\epsilon_0) 
{\d V(z,y)\over \d \nu(z)} \, d\sigma(z)\\
&&+ \int_D \Delta_z f_{x,0}(|z-x|^2+\epsilon_0) V(z,y) \, d z\\
&& - \int_D f_{x,0}(|x-z|^2+\epsilon_0) C_D \, d z.
\end{eqnarray*}
By the definition of $f_{x,0}$ and by (2.12), one can easily see 
that $I_{11}(x_0,y)$
satisfies the estimate (5.2).

 Finally, for convenience, we assume that
$n/2$ is a positive integer (otherwise we use $(n+1)/2$ instead).  Then
\begin{eqnarray*}
\lefteqn{I_{12}(x,y)}\\
&=& \int_{\d D}\biggl ({\d \Gamma(x-z)\over d\nu(z)}-
 {\d f_{x, 0}(|z-x|^2+\epsilon_0)\over \d \nu(z)}\biggr )V(z, y) \, d \sigma(z)\\
&=&  (n-2)\int_{\d D} \biggl ({\langle \nu(z), z-x\rangle \over |x-z|^n}-
{\langle \nu(z), z-x\rangle \over (|x-z|^2+\epsilon_0)^{n/2}}\biggr )
V(z,y) \, d\sigma(z)\\
&=&  (n-2)\sum_{m=0}^{n/2-1} \int_{\d D} \Big({\epsilon_0 \langle \nu(z),
 z-x\rangle 
\over |x-z|^{n-2m} (|x-z|^2+\epsilon_0)^{1+m}}
\Big) V(z,y) \,
d\sigma(z) . \\
\end{eqnarray*}
By the definition of $\epsilon_0$, one can easily 
see that $I_{12}(x_0,y)$ satisfies
the desired estimate (5.2).
Therefore, combining the above estimates,
 the proof of Theorem 5.2 is complete.\epf
\smallskip \\

As a corollary, we have the following result.

\begin{theorem}\sl Let $0<p\le 1$ and let $D$ be a bounded domain with
 $C^{\alpha(p,\epsilon)}$ boundary. Then
$$
\left|{\d^2 V(x, y)\over \d x_i \d x_j}-\sum_{|\alpha|\le n_p}
{\d^{\alpha}\over \d y^{\alpha}}{\d^2 V(x, x_0)
\over \d x_i \d x_j} (x-x_0)^{\alpha}\right|
\le C_{p,n ,\epsilon}{\delta^{n(1/p-1)+\gamma(\epsilon)}\over
|x-x_0|^{n+n(1/p-1)\gamma(\epsilon)}} \leqno (5.3)
$$
for all $x\in D\setminus 2 Q$ and $y, x_0\in Q\subset 2Q\subset D$
and any $\epsilon>0$, where $\gamma(\epsilon)=n(1/p-1)-[n(1/p-1)]+\epsilon/2
\le 1$.
\end{theorem}

In conclusion, combining Theorems 4.1 and 5.1, the proof of Theorem 1.2 is
complete. \epf

With the same argument, the Neumann problem in Theorem 1.3 can be proved
 by using Theorem 5.3 and the argument for proving Theorem 4.3. 
We leave the details for the interested reader.  The proof of Theorem
1.3 is complete.\epf

\section{Counterexamples}

In this section, we shall give some examples to show that the
hypothesis on smoothness of $\d D$ in Theorem 1.2 is essentially sharp.
We first prove the following lemma.

\begin{lemma} \sl Let $D$ be a bounded domain 
in $\RR^n$ with $C^1$  boundary. If
(integration against the kernel)
${\d^2 G\over \d x_i \d x_j}$ is bounded from $h^p_z(D)$ to $h_r^p(D)$
then, for any $x_0$ near $\partial D$ with $r(x_0)=\delta(x_0)/C(n)$, we have
$$
r(x_0)^ {pn_p-n(1-p)} \int_{D\setminus B(x_0, r(x_0))} \biggl |{\d^{2+n_p} G\over
\d x_i \d x_j \d y^{\alpha} }(x, x_0)\biggr |^p \, dx \le C_p \leqno(6.1)
$$
for all multi-indices $\alpha$ with $|\alpha|\equiv n_p=[n(1/p-1)]+1$ and $p<1/2$. 
\end{lemma} 

\proof 
Let $x_0\in D$ be near the boundary  and let $
r(x_0)=\delta(x_0)/C(n)$
with $C(n) > > 1$ a constant to be chosen so that there is
cube $Q(x_0)$ with center at $x_0$ with 
$B(x_0, r(x_0))\subset 2Q(x_0) \subset D$.
We consider the function
$$
g_{x_0}(x)= \phi_{r(x_0)}(x)= {1\over r(x_0)^n} \phi\biggl
({x-x_0\over r(x_0)}
\biggr ) . \leqno(6.2)
$$
[Here $\phi$ is a radial bump function as usual.]
For each $0<p<<1$, we define 
$$
a(x) = r(x_0)^{n_p-n(1/p-1)}{\d^ {n_p} g_{x_0}(x)\over \d x^{\alpha}}
. \leqno(6.3)
$$
 It is easy to show that $a$ is a $p$-atom with
support in $B(x_0, r(x_0))\subset 2Q\subset D$.
Thus, for all $x\in D\setminus B(x_0, r(x_0))$, since $G(x, \cdot)$ is
harmonic in $B(x_0, r(x_0))$, we have
\begin{eqnarray*}
\lefteqn{\int_D {\d^2 G(x,y) \over
\d x_i \d x_j} a(y) dy}\\
& =& \pm  r(x_0)^{n_p-n(1/p-1)} \int_D {\d^{2+n_p} G(x,y) \over
\d x_i \d x_j \d y^{\alpha} } g_{x_0}(y) dy\\
&=& \pm  r(x_0)^{ n_p-n(1/p-1)}{\d^{2+n_p} G(x,y) \over
\d x_i \d x_j \d y^{\alpha} }\biggr |_{y=x_0}\\
\end{eqnarray*}
for all multi-indices $\alpha$ with $|\alpha|=n_p$. Therefore
the fact ${\d^2 G\over \d x_i \d x_j}$ maps $h^p_z(D)$ to $h_r^p(D)$
boundedly
 implies that (6.1) holds, and the proof of Lemma 6.1 is complete.\epf

Let $B_2$ be the unit disc in $\RR^2$ and let $\psi$ 
be a conformal map from $B_2$ to some domain $D \subseteq \CC \approx \RR^2$..
We write $D=\psi(B_2)$, and $\varphi(z)=\psi^{-1}(z)$. Then the Green's function
for $-\Delta$ on $D$ is
$$
G_D(z,w)={1\over 2\pi} \log {|1-\varphi(z) \overline{\varphi}(w)|\over |\varphi(z)-\varphi(w)|}
$$
for all $z, w\in D$.  Now we have

\begin{proposition} \sl With notation above,
$$
4\pi {\d^{2+2_p} G(z,w) \over  \d w^{2_p} \d z^2 }=
-\varphi''(z){\d^{2_p} \over \d w^{2_p}}
\biggl ({1 \over \varphi(z)-\varphi(w)}\biggr ) +\varphi'(z)^2
{\d^{2_p} \over \d w^{2_p}} \biggl ({1\over
(\varphi(z)-\varphi(w))^2}\biggr )
$$
for all $z\ne w$, where $2_p=[2(1/p-1)]+1$.
\end{proposition}
\proof  We calculate that
\begin{eqnarray*}
\lefteqn{4\pi {\d^{2+2_p} G(z,w) \over \d w^{2_p} \d z^2 }}\\
&=& {\d^{2+2_p} \over \d z^2 \d w^{2_p}}(\log |1-\varphi(z) \overline{\varphi}(w)|^2 -\log|\varphi(z)-\varphi(w)|^2)\\
&=&{\d^{2_p} \over \d w^{2_p}}
\Big({-\varphi''(z) \overline{\varphi(w)}\over 1-\varphi(z) \overline{\varphi}(w)}-{\varphi'(z)^2 \overline{\varphi(w)}^2
 \over (1-\varphi(z) \overline{\varphi}(w))^2}\\
&& -{\varphi''(z)  \over \varphi(z)-\varphi(w)} +{\varphi'(z)^2 
\over (\varphi(z)-\varphi(w))^2}\biggr )\\
&=& -\varphi''(z){\d^{2_p} \over \d w^{2_p}}
\biggl ({1 \over \varphi(z)-\varphi(w)}\biggr ) +\varphi'(z)^2 
{\d^{2_p} \over \d w^{2_p}} \biggl ({1\over
(\varphi(z)-\varphi(w))^2}\biggr )
\end{eqnarray*}
and the proof is complete. \epf

Without loss of generality, we may assume henceforth that $2(1/p-1)>2$ is an
integer.

\begin{proposition}\sl Suppose that
 $\psi \in C^{\beta}(\overline{B_2})$ (the same for $\varphi$) with
$\beta=2(1/p-1)-\epsilon$ (for this last, it is sufficient that
$\partial D$ be $C^{2(1/p - 1) - \epsilon/2}$. 
Assume that ${\d^2 G\over \d x_i \d x_j}$ induces
a bounded operator from $h_z^p(D)$ to $h_r^p(D)$. Then
$$
r(x_0)^p|\varphi^{(2_p)}(x_0)|^p C_{\varphi, p, x_0}
\le C_p(1+ \|\varphi\|_{\Lambda_{\beta}})^{2_p},
$$
where
$$
C_{\varphi, p, x_0} \equiv
 \int_{D\setminus B(x_0, r(x_0))}\Big|{ 2\varphi'(z)^2 \over
(\varphi(z)-\varphi(x_0))^3 }-{\varphi''(z) \over (\varphi(z)-\varphi(x_0))}\Big|^p \, dA(z)  .
\eqno (\dagger)
$$
\end{proposition}
\proof Now
\begin{eqnarray*}
\lefteqn{\biggl | 4\pi {\d^{2+2_p} G(z,w) \over \d z^2 \d w^{2_p}}\biggr | }\\
&\ge& \biggl |{2\varphi'(z)^2 \varphi^{(2_p)}(w)\over (\varphi(z) -\varphi(w))^3}
-{\varphi''(z) \varphi^{(2_p)}(w)\over
(\varphi(z)-\varphi(w))^2}\biggr |\\
&& - \sum_{k=1}^{2_p-2} {C_p (1+|\varphi|_{C^{2_p-2}(\Dbar)})^{2_p}
\over |\varphi(z)-\varphi(w)|^{3+ 2_p-k}}
- {C_p (1+\|\varphi\|_{\Lambda_{\beta}(\Dbar)})\delta(w)^{-\epsilon}
\over |\varphi(z)-\varphi(w)|^{4}}\\\
&\ge& \biggl |{\varphi^{(2_p)}(w)\over
|\varphi(z)-\varphi(w)|^2}\biggr | 
\biggl |
{2\varphi'(z)^2 \over (\varphi(z) -\varphi(w))} 
-\varphi''(z) \biggr |\\
&&- {C_p \over |\varphi(z)-\varphi(w)|^{2+2_p}}
-{C_p \over  \delta(w)^{\epsilon} |\varphi(z)-\varphi(w)|^4}.
\end{eqnarray*}

We conclude that
\begin{eqnarray*}
\lefteqn{\biggl |{\varphi^{(2_p)}(w)\over
|\varphi(z)-\varphi(w)|^2}\biggr | 
\biggl |
{2\varphi'(z)^2 \over (\varphi(z) -\varphi(w))} 
-\varphi''(z) \biggr |}\\
&\le&  4\pi \biggl | {\d^{2+2_p} G(z,w) \over \d z^2 \d w^{2_p}}\biggr |
+{C_p\over |\varphi(z)-\varphi(w)|^{2+2_p}}+
{C_p \|\varphi\|_{\Lambda_{\beta}} \over 
\delta(w)^{\epsilon} |\varphi(z)-\varphi(w)|^4}.
\end{eqnarray*}
Also 
$$
r(x_0)^{p}
\int_{D\setminus B(x_0, r(x_0))}
 {1\over |\varphi(z)-\varphi(x_0)|^{p(2+n_p)}} \, dA(z)
\le {C_p r(x_0)^{p} \over
 r(x_0)^{2+p -2}}
= C_p
$$
and 
$$
r(x_0)^{p}
\int_{D\setminus B(x_0, r(x_0))} {1 \over \delta(z)^{p \epsilon}
|\varphi(z)-\varphi(x_0)|^{4p}} \, dA(z)
\le  C_p r(x_0)^{p-p\epsilon} 
\le  C_p.
$$
Since $2(1/p-1)$ is integer, we have
$$
p 2_p-2(1-p)= p(2(1/p-1)+1)-2(1-p)=p.
$$
Combining this and all those estimates with Lemma 6.1 and Propositions 6.2, the proof of the proposition is complete. \epf

\begin{proposition} \sl For any $\epsilon>0$, there is a bounded
domain $D$ in $\RR^2$ with $C^{2(1/p-1)-\epsilon}$ boundary and 
the operator induced by ${\d^2 G\over \d x _i \d x_j}$
is not bounded from $h^p_z(D)$ to $h^p_r(D)$.
\end{proposition}

\proof Let 
$$
\psi(z)=z+ \eta (1-z)^{2(1/p-1)-\epsilon}
$$
with $0<\eta=\eta(p) <<1$ sufficiently small, depending on $p$. Clearly, 
$$
\psi'(z) =1+ \eta (2(1/p-1)-\epsilon)(1-z)^{2/p-3-\epsilon}
$$
It is easy to see that $\psi$ is a conformal map from $B_2$ 
onto $\psi(B_2)$ provided that $0<\eta=\eta(p)$ is sufficiently small
and $\psi
\in \Lambda_{\beta}(B_2)$ where $\beta=2(1/p-1)-\epsilon$. Let
$D=\psi(B_2)$ and $\varphi(w)= \psi^{-1}(w)$.  It is clear that
$$
C_{\varphi, p, \psi(1-\delta)}>1/C>0
$$
(refer to equation $(\dagger)$).
Now, since $2_p=2(1/p-1)+1$, we have
$$
\psi^{(2_p)}(z)=c_{p ,\epsilon} (1-z)^{-\epsilon-1}
$$
where $c_{p,\epsilon} \ne 0$. 
[Note that superscripts in parentheses denote derivatives.]  
We let $y_0=1-\delta$. Then
$$
\delta ^p \biggl |\psi^{(2_p)}(y_0)\biggr |^p=|c_{p,\epsilon}|^p\delta^p 
\delta^{-\epsilon p-p}
=|c_{p,\epsilon}|^p \delta^{-p\epsilon}
\to \infty
$$
as $\delta \to 0^+$, since $\varphi$ has boundary
behavior similar to that of $\psi$.
  Thus
$$
\delta^{-p\epsilon} \le {C_p\over |c_{p,\epsilon}|^p}
\delta^p \biggl |\varphi^{(2_p)}(\psi(1-\delta)) \biggr | .
$$

Seeking a contradiction, we suppose that  
${\d^2 G\over \d x_i \d x_j}$ induces a bounded 
operator from $h_z^p(D)$ to $h^p_r(D)$.
Then, by Propositions 6.2 and 6.3, we have
$$
\delta^{-p\epsilon} \le {C_p\over |c_{p,\epsilon}|^p}
\delta^p |\varphi^{(2_p)}(\psi(1-\delta))|
< C_{p,\epsilon}(1+\|\varphi\|_{\Lambda_{\beta}})^{2_p}<\infty
$$
for any $0<\delta<<1$.
 This is a contradiction as $\delta\to 0^+$.
Thus the operator induced by ${\d^2 G\over \d x_i \d x_j}$ is 
not bounded from $h_z^p(D)$ to $h^p_r(D)$,
and the proof is complete. \epf

\end{document}